\documentclass[9pt,amsfonts]{article}
\usepackage{graphicx,amssymb,eucal,mathrsfs}
\usepackage{latexsym,amsmath,amscd,xcolor}
\usepackage[all]{xy}
\usepackage{ntheorem,tikz}
\usepackage{pgf,etex,pstricks}
\usepackage{enumitem}
\usepackage[hypertexnames=false,backref=page,pdftex,
 	pdfpagemode=UseNone,
 	breaklinks=true,
 	extension=pdf,
 	colorlinks=true,
 	linkcolor=blue,
 	citecolor=red,
 	urlcolor=blue,
 ]{hyperref}

\newcommand{\Z}{\mathbb{Z}}

\newcommand{\out}{{\rm Out}}
\newcommand{\aut}{{\rm Aut}}
\newcommand{\inn}{{\rm Inn}}
\newcommand{\hp}{\ensuremath{\mathbb H}^2}

\newcommand{\diff}{{\rm Diff}}
\newcommand{\mcg}{{\rm MCG}}

\newtheorem{main}{Theorem}

\newtheorem{question}[main]{Question}
\newtheorem{prop}[main]{Proposition}
\newtheorem{lemma}[main]{Lemma}
\newtheorem{rem}[main]{Remark}
\newtheorem{cor}[main]{Corollary}

\newtheorem{claim}[main]{Claim}
\newtheorem{imain}{Theorem}

\theoremstyle{definition}

\newtheorem{defi}[main]{Definition}

\title{On finitely generated normal subgroups of K\"ahler groups}
\date{July 2021}
\author{Francisco Nicol\'as}

\begin{document}

\maketitle

\begin{abstract} 
\noindent We prove that if a surface group embeds as a normal subgroup in a K\"ahler group and the conjugation action 
of the K\"ahler group on the surface group preserves the conjugacy class of a non-trivial element, then the K\"ahler group is 
virtually given by a direct product, where one factor is a surface group. Moreover we prove that if a one-ended hyperbolic group 
with infinite outer automorphism group embeds as a normal subgroup in a K\"ahler group then it is virtually a surface group. 
More generally we give restrictions on normal subgroups of K\"ahler groups which are amalgamated products or HNN extensions.
\end{abstract}


\section{Introduction}

\noindent Let $G$ be a finitely generated group. Bass-Serre Theory (see~\cite{SW} and~\cite{Serre}) establishes a dictionary 
between decompositions of $G$ as an amalgamated product or an HNN extension and actions 
of $G$ on a simplicial tree without inversions which are transitive on the set of edges. Of course, the theory 
also deals with more complicated graphs of groups but we will mainly deal with amalgamated products and HNN extensions.
If $G$ acts on a tree $T$ and we have finitely generated groups $\Gamma$ and $Q$ that fit in a short exact sequence 
\begin{gather}\label{secp}
 \xymatrix{
 1\ar[r] & G \ar[r] &\Gamma \ar[r] &Q \ar[r] &1,
 }
\end{gather}
one can ask whether the action of $G$ on $T$ can be extended to $\Gamma$. 
If the center of $G$ acts trivially on $T$, we obtain an 
induced action of the group $\inn(G)$ of inner automorphisms of $G$ on $T$. As before, one can ask whether this induced action can 
be extended to a larger group of automorphisms of $G$. Of course both questions are related, since for groups $\Gamma$ and 
$G$ as in \eqref{secp} the conjugation action of $\Gamma$ on $G$ induces a map $\Gamma\rightarrow \aut(G)$. Several authors 
have studied the latter question. Karrass, Pietrowski and Solitar~\cite{KPS} studied the case of an amalgamated product $A\ast_CB$, 
where $C$ is maximal among all its conjugates in $A$ and $B$. Pettet~\cite{Pettet} studied the general situation of a graph of groups with 
a more restrictive condition on the edge stabilizers of its Bass-Serre tree (\textit{``edge group incomparability hypothesis''}) which 
is equivalent to the \textit{conjugate maximal condition} 
of Karrass, Pietrowski and Solitar when there is one orbit of edges. A particular case of this situation was studied by Gilbert, 
Howie, Metaftsis and Raptis~\cite{GHMR} where they proved that the action of the Baumslag-Solitar group on its Bass-Serre tree can be 
extended to its whole group of automorphisms in many cases.
For a one-ended torsion-free hyperbolic group 
$G$ which is not virtually a surface group, Sela~\cite{Sela} proved the existence of a ``canonical'' tree $T$ on which $G$ acts.
By studying the action of $\aut(G)$ on $\partial G$, Bowditch~\cite{Bowditch} proved that the action of $G$ on $T$ can 
be extended to $\aut(G)$.

We will apply these results to study finitely generated normal subgroups of K\"ahler groups. 
Recall that a group is called a \textit{K\"ahler group} if it can be realized as the fundamental group of a compact K\"ahler manifold, 
a classical reference on this subject is~\cite{ABCKT}. Some examples of K\"ahler groups are Abelian groups of even rank, 
surface groups (groups that can be realized as the fundamental group of a closed hyperbolic surface) and uniform lattices of $PU(n,1)$ 
(the group of holomorphic isometries of complex hyperbolic $n$-space). 
The main ingredient to apply these results about actions on trees that extend to a group of automorphisms is a classical result of 
Gromov and Schoen~\cite{GS} about K\"ahler groups acting on trees (see also \cite{Sun}).

First, we will study short exact sequences as in \eqref{secp}, where $\Gamma$ is a  K\"ahler group and $G$ is a surface group, 
\textit{i.e.}, when $G$ can be realized as the fundamental group of a closed hyperbolic surface $S$.
The first main result of this work is the folllowing:

\begin{imain}\label{t}
 Suppose that the conjugation action of $\Gamma$ on $\pi_1(S)$ preserves the conjugacy class of a simple 
 closed curve in $S$. Then there is a finite index subgroup $\Gamma_1$ of $\Gamma$ containing $\pi_1(S)$
 such that the restricted short exact sequence
 \begin{gather*}
  \xymatrix{
  1\ar[r] &\pi_1(S) \ar[r] &\Gamma_1 \ar[r] &Q_1 \ar[r] &1
  }
 \end{gather*}
 splits as a direct product (where $Q_1$ is the image of $\Gamma_1$ in $Q$).
\end{imain}
A particular case arises when a K\"ahler extension by a surface group has an Abelian quotient. In that case, using classical results on 
the mapping class group we have that for such a short exact sequence, the hypothesis of Theorem \ref{t} is 
virtually fulfilled if the rank of the image of $Q\rightarrow\out(\pi_1(S))$ is greater than $1$. From this we can deduce that a map
$\psi:Q\rightarrow\out(\pi_1(S))$ with $Q$ any Abelian group, is the monodromy map of a short exact sequence as in \eqref{secp}, with 
$\Gamma$ a K\"ahler group and $G=\pi_1(S)$, only if it has finite image. This result is due to Bregman and Zhang (see~\cite{BZ}) who 
established it using different techniques. We recover it as a corollary of Theorem \ref{t} (see Section \ref{qab}). 

Going back to the case of an arbitrary quotient $Q$, Theorem \ref{t} says that
if the image of the monodromy map $Q \rightarrow\out(\pi_1(S))$ is 
infinite, it cannot preserve the isotopy class of any simple closed curve in $S$.
It is well known that this conclusion holds when the short exact sequence is
induced by a surjective holomorphic submersion (which is not a locally trivial holomorphic fiber bundle) between compact 
complex manifolds, whose fibers are homeomorphic to the closed hyperbolic surface $S$ and whose base is either projective or a 
Riemann surface of finite type. This is a theorem due to Shiga~\cite{Shiga}. Hence, our theorem can be seen as a generalization 
of Shiga's result (in the case of compact bases).

We now go back to the study of extensions as in \eqref{secp} where $G$ need not be a surface group. 
The folllowing three results are an application of our study on K\"ahler groups and the works~\cite{Bowditch, GHMR, KPS, Pettet, Sela}. 
Recall that a group $H$ is indecomposable if for any decomposition $H=H_1 \ast H_2$ as a free product, $H_1$ or $H_2$ is trivial.

\begin{imain}\label{kgfp}
 Let $G$ be a group that splits as a non-trivial free product $A\ast B$ with $A$ and $B$ indecomposable and not infinite cyclic. 
 Then $G$ does not embed as a normal subgroup in a K\"ahler group.
\end{imain}

\begin{imain}\label{kgbs}
 Let $G$ be the Baumslag-Solitar group
 $$\langle x,t | tx^{p}t^{-1}=x^{q}\rangle,$$
 where $p, q$ are integers with $p,q > 1$ and such that neither one is a multiple of the other.
 Then $G$ does not embed as a normal subgroup in a K\"ahler group.
\end{imain}

Finally, our last result deals with the case where the normal subgroup $G$ is a one-ended hyperbolic group.

\begin{imain}\label{kghg}
 Suppose that $G$ is a one-ended torsion-free hyperbolic group that embeds as a normal subgroup in a K\"ahler group $\Gamma$.
 If $\out(G)$ is infinite, then $G$ is virtually a surface group. 
\end{imain}
We observe that in this last theorem there is no hypothesis on the monodromy map $\Gamma\rightarrow \out(G)$. We only need to 
know that $\out(G)$ is infinite.

Theorem \ref{kgfp} raises the following natural question:

\begin{question}\label{freeproductkahler}
  Is there a K\"ahler group admitting as a normal subgroup a non-Abelian free group or a free product $ A_1 \ast \cdots \ast A_n $ with 
  $n\geq 3$ and $A_1, \ldots ,A_n$ indecomposable and not infinite cyclic groups ?
\end{question}

The structure of this work is the following. In Section \ref{GGT},  we explain how to extend the action of a group on a tree to a 
subgroup of its group of automorphisms. In particular we study the case of surface groups.
In Section \ref{AKG}, we apply these results to the study of K\"ahler groups admitting as a normal subgroup a group acting on a tree.
Theorems \ref{t}, \ref{kgfp}, \ref{kgbs} and \ref{kghg} are proved in Section \ref{App}. Finally, in Section \ref{s4} we study the 
monodromy map of a K\"ahler extension by a surface group and establish variations on Theorem \ref{t}.


\section{Group actions on trees}\label{GGT}

\subsection{Bass-Serre dictionary}\label{BSTD}

{\bf Notations.} Suppose that $G$ acts on a tree $T$ and let $e$ and $v$ be an edge and a vertex of $T$ respectively. We will denote by 
$G_e$ and $G_v$ the subgroups of $G$ given by the respective stabilizers of $e$ and $v$.

\noindent The dictionary established by Bass-Serre theory is given by the following results (see~\cite{SW, Serre}):

\begin{main}[Serre]
 Let $G$ be a group acting on a tree $T$ without inversions in such a way that the action is transitive on the set of edges. 
 Let $e=(v_1,v_2)$ be an edge of $T$.
 \begin{enumerate}
  \item If $v_1$ and $v_2$ are in different orbits, then $G$ splits as the amalgamated product $G_{v_1}\ast_{G_e}G_{v_2}$.
  \item If $v_1$ and $v_2$ are in the same orbit, then $G$ splits as the HNN extension $G_{v_1}\ast_{G_e,\theta}$ where 
  $\theta: G_e \rightarrow G_{v_1}$ is the monomorphism given by the conjugation by an element $t$ in $G\setminus G_{v_1}$ 
  that sends $v_1$ to $v_2$.
 \end{enumerate}
\end{main}

\begin{main}[Serre]
 Let $G$ be a group that splits as an amalgamated product or as an HNN extension. Then, $G$ acts on a tree $T$ without inversions 
 and this action is transitive on the set of edges.
 \begin{enumerate}
  \item If $G$ splits as the amalgamated product $A\ast_CB$, then the set of vertices of $T$ is the 
  disjoint union of left cosets $G/A\sqcup G/B$ and the set of edges of $T$ is the set of left cosets $G/C$. 
  The adjacency is given by the maps $G/C\rightarrow G/A$ and $G/C\rightarrow G/B$.
  \item If $G$ splits as the HNN extension $A\ast_{C,\theta}$, where $\theta: C\rightarrow A$ is the monomorphism given by the 
  conjugation by an element $t$ in $G\setminus A$, then the set of vertices of $T$ is the set of 
  left cosets $G/A$ and the set of edges of $T$ is the set of left cosets $G/C$. 
  The adjacency is given by the maps $\iota_1:G/C\rightarrow G/A$ and $\iota_2:G/C\rightarrow G/A$ that send $xC$ to $xA$ 
  and $xC$ to $xt^{-1}A$ respectively.
 \end{enumerate}
\end{main}
In both cases we will say that \textit{$G$ splits over the subgroup $C$}.
The tree associated in this way to a group that splits over a subgroup $C$
is called the Bass-Serre tree of $G$. The action of $G$ on its Bass-Serre tree $T$ has the following properties:
\begin{enumerate}
 \item It is minimal, \textit{i.e.} there is no $G$-invariant proper subtree in the Bass-Serre tree $T$.
 \item If $G$ splits as an amalgamated product $A \ast_C B$ and $C$ is properly contained in $A$ and $B$, or 
 if $G$ splits as an HNN extension $A\ast_{C,\theta}$, and $C$ and $\theta(C)$ are properly contained in $A$
 then the action on the boundary of the tree has no fixed points. Recall that the 
 boundary of a tree is given by the set of infinite paths without backtracking starting at a fixed point of the tree 
 (this is one definition of the boundary of a tree among many others).
 \item If $G$ splits as an amalgamated product $A\ast_CB$ such that $[A:C]>2$ or $[B:C]>2$, or if $G$ splits as an HNN extension 
 $A\ast_{C,\theta}$ such that $[A:C]>1$ or $[A:\theta(C)]>1$, then the Bass-Serre tree of $G$ is not a line.
\end{enumerate}


\subsection{Extending the action on the Bass-Serre tree}\label{EA}

Let $G$ be a finitely generated group that splits as an amalgamated product or an HNN extension and let $T$ be its Bass-Serre tree. 
If the center $Z(G)$ of $G$ acts 
trivially on $T$, the action of $G$ on $T$ factors through the quotient $G/Z(G)$. Since this quotient is isomorphic to the group 
$\inn(G)$ of inner automorphisms of $G$, we obtain an induced action of $\inn(G)$ on $T$. 
The aim of this section is to extend this
induced action of $\inn(G)$ on $T$ to a larger group of automorphisms of $G$.

\begin{defi}
 We define the subgroup $\aut_T(G)$ of $\aut(G)$ by the following property:
 an automorphism $\varphi:G\rightarrow G$ is an element of $\aut_T(G)$ if for any edge $e=(v_1,v_2)$ of $T$ there is an 
 element $x$ in $G$ such that 
 $$\varphi(G_{v_1})=x G_{v_1} x^{-1},\quad \varphi(G_{v_2})=x G_{v_2} x^{-1}, \quad \mbox{and} \quad \varphi(G_{e})=x G_{e} x^{-1}.$$
\end{defi}

\begin{lemma}\label{iso}
 Suppose that the stabilizer of any vertex of $T$ is its own normalizer and let $\varphi:G\rightarrow G$ be an element of 
 $\aut_T(G)$. Then, $\varphi$ induces an isometry $\overline{\varphi}:T\rightarrow T$ such that for all $g$ in $G$
 \begin{equation}\label{gc}
  \overline{\varphi}(g\cdot \underline{\,\,\,}\,)=\varphi(g) \cdot \overline{\varphi}(\,\underline{\,\,\,}\,). 
 \end{equation}
 Moreover, if $\varphi$ is the inner automorphism given by the conjugation by $x$, then $\overline{\varphi}$ is given by 
 the action of $x$ on $T$.
\end{lemma}

\begin{rem}\label{remlemiso}
 The hypothesis of Lemma \ref{iso} on the vertex stabilizers implies that the center of $G$ acts trivially on $T$.
\end{rem}

\noindent {\it Proof of Lemma \ref{iso}.} 
 Let us fix an edge $e=(v_1,v_2)$. We warn the reader that the following construction of $\overline{\varphi}$ depends on the choice of $e$. 
 We write
 $$\varphi(G_{v_1})=x G_{v_1} x^{-1}, \quad \varphi(G_{v_2})=x G_{v_2} x^{-1}, \quad \mbox{and} \quad \varphi(G_{e})=x G_{e} x^{-1}.$$
 We define an isometry $\overline{\varphi}:T\rightarrow T$ as follows. If $v$ is a vertex of $T$ we write $v=g\cdot v_i$ for some 
 index $i$ and define
 $$\overline{\varphi}(v)=\varphi(g) x\cdot v_i.$$
 If $f$ is an edge of $T$ we write $f=g\cdot e$ and define
 $$\overline{\varphi}(f)=\varphi(g) x\cdot e.$$
These definitions are independent of the choice of the element $g$. In the case of an HNN extension, we must check that the action on 
a vertex $v$ is independent of whether we represent $v$ as the image of $v_1$ or $v_2$ (as $v_1$ and $v_2$ are in the same orbit). 
But if $v_2=t\cdot v_1$ one checks that $t^{-1}x^{-1}\varphi (t)x$ normalizes $G_{v_1}$, and thus by the 
hypothesis of the lemma we get that $t^{-1}x^{-1}\varphi (t)x$ lies in $G_{v_1}$. Using this observation one proves that 
$\overline{\varphi}(v)$ is well defined. In this way we obtain two bijections, one of the set of vertices of $T$ and one of the set 
of edges of $T$, which define an isometry of $T$.
If there is another element $y$ in $G$ such that 
$$\varphi(G_{v_1})=y G_{v_1} y^{-1}, \quad \varphi(G_{v_2})=y G_{v_2} y^{-1}, \quad \mbox{and} \quad \varphi(G_{e})=y G_{e} y^{-1}$$
then $x^{-1}y$ is in the normalizer of $G_{v_1}$ and $G_{v_2}$. 
Since $G_{v_1}$ and $G_{v_2}$ are their own normalizers we get that $x^{-1}y$ is in the intersection of  
$G_{v_1}$ and $G_{v_2}$, which is $G_e$. This implies that $\overline{\varphi}$ is well defined and it has the desired property
$$\overline{\varphi}(g\cdot \underline{\,\,\,}\,)=\varphi(g) \cdot \overline{\varphi}(\,\underline{\,\,\,}\,).$$
Finally if $\varphi$ is the inner automorphism given by the conjugation by $x$ we get that 
$$\overline{\varphi}(g\cdot e)=\varphi(g) x \cdot e = (xgx^{-1})x \cdot e = x\cdot (g \cdot e).$$
Similarly, for $i=1,2$ we obtain that $\overline{\varphi}(g\cdot v_i)= x\cdot (g \cdot v_i)$, which implies that $\overline{\varphi}$ 
is given by the action of $x$ on $T$.
\hfill $\Box$

\begin{defi}
An isometry $\overline{\varphi}:T\rightarrow T$ induced by an automorphism $\varphi$ of $G$ that satisfies \eqref{gc} is called 
\textit{$G$-compatible}. An extension of the induced  action of $\inn(G)$ on $T$ to a subgroup $\Lambda$ of $\aut(G)$ is called 
\textit{$G$-compatible} if every automorphism of $\Lambda$ defines a $G$-compatible isometry of $T$.
\end{defi}

\begin{lemma}\label{giso}
 Suppose that the stabilizer of any vertex of $T$ is its own normalizer. Then the induced action of $\inn(G)$ on $T$ extends to a 
 $G$-compatible action of $\aut_T(G)$ on $T$.
\end{lemma}

\noindent {\it Proof.} 
 Let us fix an edge $e=(v_1,v_2)$.
 By Lemma \ref{iso}, it suffices to show that for all $\varphi,\psi$ in $\aut_T(G)$ 
 the isometries $\overline{\varphi}$ and $\overline{\psi}$ defined as before satisfy:
 $$\overline{\varphi^{-1}}=\overline{\varphi}^{-1} \quad \mbox{and} \quad 
 \overline{\varphi\circ\psi}= \overline{\varphi}\circ\overline{\psi}.$$
 
 \noindent Let $\varphi$ be an element in $\aut_T(G)$ such that 
 $$\varphi(G_{v_1})=x G_{v_1} x^{-1},\quad \varphi(G_{v_2})=x G_{v_2} x^{-1}, \quad \mbox{and} \quad \varphi(G_{e})=x G_{e} x^{-1}.$$
 Then, if $y=\varphi^{-1}(x^{-1})$ we get that
 $$\varphi^{-1}(G_{v_1})=y G_{v_1} y^{-1}, \quad \varphi^{-1}(G_{v_2})=y G_{v_2} y^{-1}, 
 \quad \mbox{and} \quad \varphi^{-1}(G_{e})=y G_{e} y^{-1}.$$
 Therefore the isometry $\overline{\varphi^{-1}}$ is given by 
 $$\begin{array}{lll}
   \overline{\varphi^{-1}}(g \cdot v_i)&=\varphi^{-1}(g)y\cdot v_i & \mbox{for $i=1,2$};\\
   \overline{\varphi^{-1}}(g \cdot e)&=\varphi^{-1}(g)y\cdot e. &
 \end{array}$$ 
 A direct computation shows that $\overline{\varphi^{-1}}=\overline{\varphi}^{-1}$.

\noindent Now, let $\varphi$ and $\psi$ be two elements of $\aut(G)$ such that 
$$\varphi(G_{v_1})=x G_{v_1} x^{-1},\quad \varphi(G_{v_2})=x G_{v_2} x^{-1}, \quad \mbox{and} \quad \varphi(G_{e})=x G_{e} x^{-1},$$
and 
$$\psi(G_{v_1})=y G_{v_1} y^{-1},\quad \psi(G_{v_2})=y G_{v_2} y^{-1}, \quad \mbox{and} \quad \psi(G_{e})=y G_{e} y^{-1}.$$
Then, if $z=\psi(x)y$ we get that 
$$\psi\circ\varphi(G_{v_1})=z G_{v_1} z^{-1},\quad \psi\circ\varphi(G_{v_2})=z G_{v_2} z^{-1}, 
\quad \mbox{and} \quad \psi\circ\varphi(G_{e})=z G_{e} z^{-1}.$$
Hence the isometry $\overline{\psi\circ\varphi}$ is given by 
 $$\begin{array}{lll}
   \overline{\psi\circ\varphi}(g \cdot v_i)&=\psi\circ\varphi(g)z\cdot v_i & \mbox{for $i=1,2$};\\
  \overline{\psi\circ\varphi}(g \cdot e)&=\psi\circ\varphi(g)z\cdot e. &
 \end{array}$$ 
Once again, a direct computation shows that $\overline{\psi\circ\varphi}=\overline{\psi}\circ\overline{\varphi}$.
\hfill $\Box$

Note that the isometry $\overline{\varphi}:T\rightarrow T$ defined in Lemma \ref{iso} depends on the selected edge $e=(v_1,v_2)$ when 
$\varphi$ is not an inner automorphism of $G$. Therefore, the action of $\aut_T(G)$ on $T$ also depends on this edge.

\begin{defi}\label{defini}
 We denote by $\Lambda_T$ the subgroup of $\aut(G)$ that preserves the conjugacy class of each vertex stabilizer and each edge 
 stabilizer of the Bass-Serre tree $T$ of $G$, \textit{i.e.} an automorphism $\varphi$ of $G$ is an element of $\Lambda_T$ if for any 
 edge $e=(v_1,v_2)$  there are elements $x,y,z$ in $G$ such that 
 $$\varphi(G_{v_1})=x G_{v_1} x^{-1},\quad \varphi(G_{v_2})=y G_{v_2} y^{-1}, \quad \mbox{and} \quad \varphi(G_{e})=z G_{e} z^{-1}.$$
\end{defi}

It is clear that $\aut_T(G)$ is always contained in $\Lambda_T$. In the following, we will see some situations where both groups 
almost coincide.


\subsection{Surface groups}\label{SG}

\noindent Let $\gamma$ be a simple closed curve in a closed hyperbolic surface $S$. If we denote by $C$ the cyclic subgroup of 
$\pi_1(S)$ generated by the homotopy class of $\gamma$, we have that  $\pi_1(S)$ splits over the cyclic subgroup $C$ as follows. 
Let us fix a base point $x_0$ in $\gamma$. If $\gamma$ 
is a separating curve, we will denote by $A$ and $B$ the fundamental groups of the surfaces obtained by cutting $S$ along $\gamma$.
Similarly, if $\gamma$ is nonseparating, we will denote by  $A$ the fundamental group of the surface obtained by cutting $S$ 
along $\gamma$. All fundamental groups are based at $x_0$. In the case where $\gamma$ is not separating, we must choose 
a ``copy'' of $x_0$ in the surface obtained by cutting $S$ along $\gamma$.

If $\gamma$ is a separating curve, $\pi_1(S)$ splits as the amalgamated product $A\ast_C B$. If $\gamma$ is a nonseparating curve, 
$\pi_1(S)$ splits as the HNN extension $A\ast_{C,\theta}$, where $\theta: C\rightarrow A$ is the monomorphism given by the conjugation 
by an element $t$ in $\pi_1(S)\setminus A$. More precisely, $t$ is the homotopy class 
of a simple closed curve in $S$ that becomes a path in $S\setminus\gamma$ joining the boundary components.

\begin{rem}\label{dualtree}
 Given a simple closed curve $\gamma$ in a closed hyperbolic surface $S$ as before, there exists a ``dual tree'' associated to $\gamma$ 
 that embeds in $\hp$. This tree coincides with the Bass-Serre tree associated to the splitting of $\pi_1(S)$ over $C$ and it can be 
 constructed as follows. Let us assume that $\gamma$ is a closed geodesic in $S$ (every non-nullhomotopic simple closed curve in $S$ 
 is homotopic to a unique closed geodesic) and let $p:\hp\rightarrow S$ be the universal covering space of $S$. We denote by $L$ 
 the set of bi-infinite geodesics in $\hp$, whose images under $p$ are equal to $\gamma$. We say that two connected components of 
 $\hp\setminus L$ are related if they are separated by exactly one element of $L$. This defines a symmetric binary relation on 
 $\hp\setminus L$.The set of vertices of the dual tree is given by the connected components of $\hp\setminus L$ and the adjacency is 
 given by this symmetric binary relation, \textit{i.e.}, the set of edges of the dual tree is given by $L$. The action of $\pi_1(S)$ on $\hp$ 
 induces an action of $\pi_1(S)$ on the dual tree, which is consistent with the action of $\pi_1(S)$ on $T$.
\end{rem}

\begin{rem}
  The group of automorphisms of $\pi_1(S)$ that preserves the conjugacy class of $\gamma$ is a subgroup of index at most $2$ of
  the subgroup of $\pi_1(S)$ made of the automorphisms $\varphi:\pi_1(S)\rightarrow\pi_1(S)$ such that 
  \begin{equation}\label{conjclasspres}
   \varphi(C)=xCx^{-1}
  \end{equation}
  for some $x$ in $\pi_1(S)$. This follows from the fact that if an automorphism $\varphi$ satisfies \eqref{conjclasspres} then 
  $\varphi([\gamma])$ must be equal to either $x[\gamma]x^{-1}$ or to $x[\gamma]^{-1}x^{-1}$.
 \end{rem}

 \begin{lemma}\label{l6}
  Let $\gamma$ be a simple closed curve in $S$ and let $T$ be the Bass-Serre tree of $\pi_1(S)$  associated to $\gamma$.
  The group of automorphisms of $\pi_1(S)$ that preserves the conjugacy class of $C$ contains $\aut_T(\pi_1(S))$ as a 
  subgroup of index at most $2$. Hence $\Lambda_T$ contains $\aut_T(\pi_1(S))$ as a subgroup of index at most $2$.
 \end{lemma}
 
 \noindent {\it Proof.} 
 Let $\Gamma$ be the subgroup of $\aut(\pi_1(S))$ that preserves the conjugacy class of $C$. We have the inclusions:
 $$\aut_{T}(\pi_{1}(S))\subset \Lambda_{T}\subset \Gamma.$$
 Hence the second point of the lemma follows from the first one.
 By the Dehn-Nielsen-Baer Theorem  we have an isomorphism between the (extended) mapping class group of $S$ and 
 the outer automorphism group of $\pi_1(S)$
 \[ \mcg^{\pm}(S)=\diff^{\pm}(S)/\diff_0(S)\rightarrow \out(\pi_1(S))=\aut(\pi_1(S))/\inn(\pi_1(S))\]
 (see Dehn~\cite{Dehn} for the original proof and Nielsen~\cite{Nielsen} for the first published proof),
 which takes the class of a diffeomorphism $f:S\rightarrow S$ to the automorphism induced by $f$ on $\pi_1(S)$ (well-defined 
 up to conjugacy).  
 Let us fix an automorphism $\varphi$ of $\pi_1(S)$ contained in $\Gamma$ and 
 let $f:S\rightarrow S$ be a diffeomorphism inducing $\varphi$. 
 
 We fix a base point $x_0$ on $\gamma$ and a collar neighborhood $U_{\gamma}$ of $\gamma$. Since the curve $f \circ \gamma$ is 
 isotopic to either $\gamma$ or $\bar{\gamma}$, one can replace $f$ by a diffeomorphism $f_1$ isotopic to $f$ such that $f_1$ 
 preserves $\gamma$ globally and fixes $x_0$. Let $\diff(S,\gamma)$ be the group of diffeomorphisms of $S$ preserving $\gamma$ 
 globally, preserving each connected component of $U_{\gamma} \setminus \gamma$ and fixing $x_0$. 
 Then, a diffeomorphism contained in $\diff(S,\gamma)$ induces an automorphism of $\pi_1(S)$ that
 preserves the subgroup $C$ and the fundamental groups of the connected components of $S\setminus \gamma$. In particular, this 
 induced automorphism lies in $\aut_T(\pi_1(S))$. 
 
 The elements of $\Gamma$ which are induced by diffeomorphisms of $\diff(S,\gamma)$ form a subgroup $\Gamma_0$ of index at most $2$. 
 Hence if $\varphi \in \Gamma_0$ we obtain that $\varphi = f_* \circ I$ for some $I \in \inn(\pi_1(S))$ and for some 
 $f \in \diff(S,\gamma)$. Therefore $\varphi \in \aut_T(\pi_1(S))$. 
 \hfill $\Box$

Recall that a subgroup $H$ of a group $G$ is called \textit{malnormal} if $xHx^{-1}\cap H$ is trivial for all $x$ in $G\setminus H$.
The fundamental groups of the surfaces obtained by cutting $S$ along $\gamma$ are malnormal in $\pi_1(S)$. This implies that the 
stabilizer of any vertex of $T$ is its own normalizer. Finally, since the center of a surface group is trivial,
we obtain the following result as a consequence of Lemmas \ref{giso} and \ref{l6}.

\begin{prop}\label{psg}
 The action of $\pi_1(S)$ on $T$ extends to a subgroup of index at most $2$ of the group of automorphisms of $\pi_1(S)$ that preserves 
 the conjugacy class of $\gamma$. Moreover, this extension is $\pi_1(S)$-compatible.
\end{prop}


\subsection{Maximal families of edge stabilizers}\label{MF}

\begin{defi}
 Let $\mathscr{F}$ be a family of subgroups of $G$. We say that $\mathscr{F}$ is maximal if  for any $2$ 
 subgroups $H,K$ in $\mathscr{F}$, $K<H$ implies $K=H$.
\end{defi}

Let $G$ be a group acting on a tree $T$ and let $v$ be a vertex of $T$. We say that an edge of $T$ 
is \textit{$v$-incident} if one of its vertices is $v$. We will be interested in the following conditions:

\begin{enumerate}[label={\bfseries C\arabic*}]
\item\label{c1} The stabilizer of an edge of $T$ is properly contained in the stabilizers of its vertices.
\item\label{c2} For any vertex $v$ in $T$ the family of stabilizers of $v$-incident edges is maximal.
\end{enumerate}

\begin{rem}
 Condition \ref{c1} implies that the stabilizer of any vertex is its own normalizer and therefore the center of $G$ acts trivially 
 on $T$ (see Remark \ref{remlemiso}).
\end{rem}

\begin{main}[Karrass, Pietrowski, Solitar,~\cite{KPS}]\label{Kar}
 Let $G$ be a group that splits as an amalgamated product $A\ast_CB$ and let $T$ be its Bass-Serre tree. 
 If the action of $G$ on $T$ satisfies conditions \ref{c1} and \ref{c2}, then the induced action of $\inn(G)$ on $T$ extends to a 
 $G$-compatible action of $\Lambda_T$ on $T$.
\end{main}

\begin{main}[Pettet,~\cite{Pettet}]\label{Pet}
 Let $G$ be a group that splits as an HNN extension $A\ast_{C,\theta}$ and let $T$ be its Bass-Serre tree. 
 If the action of $G$ on $T$ satisfies conditions \ref{c1} and \ref{c2}, then the induced action of $\inn(G)$ on $T$ extends to a 
 $G$-compatible action of $\Lambda_T$ on $T$.
\end{main}

\begin{rem}
 The original statements of Theorems \ref{Kar} and \ref{Pet} only give information about the decomposition of $\Lambda_T$ as an 
 amalgamated product in the case of Theorem \ref{Kar}, and as an HNN extension in the case of Theorem \ref{Pet}, but we can deduce 
 from the proofs of these results that such decompositions are induced by an action of $\Lambda_T$ on $T$ that extends the one  
 of $G$ on $T$.
\end{rem}

In~\cite{Pettet}, Pettet gives a more general result than the one stated above. Pettet studied groups that split as the fundamental group 
of a graph of groups with a more restrictive condition on the edge stabilizers of its Bass-Serre tree (\textit{edge group incomparability 
hypothesis}), which coincides with 
the maximality of the families of $v$-incident edge stabilizers when $G$ splits as an HNN extension. 
Gilbert, Howie, Metaftsis and Raptis studied the particular situation when the edge stabilizers are cyclic. In particular they proved the 
following result on Baumslag-Solitar groups:

\begin{main}[Gilbert, Howie, Metaftsis, Raptis]\label{bsgaut}
 Let $G$ be the Baumslag-Solitar group
 $$\langle x,t | tx^{p}t^{-1}=x^{q}\rangle,$$
 where $p, q$ are integers with $p,q > 1$ and such that neither one is a multiple of the other and let $T$ be the Bass-Serre tree of $G$. 
 Then the induced action of $\inn(G)$ on $T$ extends to a $G$-compatible action of $\aut(G)$ on $T$.
\end{main}

Note that if $G$ splits as a free product $A\ast B$,
then for any vertex $v$ in $T$ the family of stabilizers of $v$-incident edges is maximal. 
Moreover, if $A$ and $B$ are indecomposable, not infinite cyclic and not isomorphic to each other, 
then $\Lambda_T$ coincides with $\aut(G)$ (see~\cite[p. 152]{SW}). As a consequence of Theorem \ref{Kar}, we obtain the following 
result:

\begin{main}\label{fpt}
 Let $G$ be a group that splits as a free product $A\ast B$ with $A$ and $B$ indecomposable, not infinite cyclic and not isomorphic 
 to each other and let $T$ be the Bass-Serre tree of $G$. Then, the induced action of $\inn(G)$ on $T$ extends to a $G$-compatible 
 action of $\aut(G)$ on $T$.
\end{main}

 \begin{rem}\label{fpfr}
  If $G=A\ast A$, the same argument proves that $\aut(G)$ contains a subgroup of index $2$ which acts on the Bass-Serre tree of $G$. 
  The group $\aut(G)$ itself acts on the first barycentric division of the Bass-Serre tree of $G$.
 \end{rem}

  If the action of a finitely generated group $G$ on a tree $T$ induces a decomposition of $G$ as a non-Abelian free group or as a free 
  product of three (or more) indecomposable subgroups, then the study of extensions of the induced action of $\inn(G)$ on $T$ to larger 
  groups of automorphisms of $G$ could lead to an answer to Question \ref{freeproductkahler}.
 For some important contributions to the study of the group of automorphisms of a free product of a finite number of indecomposable 
 groups see for instance the work of Foux-Rabinovitch~\cite{FR1,FR2}.


\subsection{One-ended hyperbolic groups}\label{OE}

Interesting examples of groups with a non-trivial splitting are given by one-ended hyperbolic groups whose 
group of outer automorphisms is infinite. In Section \ref{SG} we saw that a surface group splits over 
a cyclic subgroup. The last assertion holds as well for a group which is virtually a surface group. When the group is a one-ended 
hyperbolic group which is not virtually a surface group, the theory 
of JSJ decompositions guarantees the existence of a non-trivial splitting over a virtually cyclic subgroup under the assumption that its 
group of outer automorphisms is infinite (see for instance~\cite[Theorem~1.4]{Levitt}). The notion of JSJ decomposition
was introduced in~\cite{Sela} by Sela for one-ended hyperbolic groups without torsion. In~\cite{Bowditch}, Bowditch 
constructed a slightly different JSJ decomposition for one-ended hyperbolic groups (not 
necessarily torsion-free) by studying cut points on the boundary of the group. Both decompositions are ``canonical'' and they coincide 
when the group is torsion-free. In the rest of this section, $G$ will denote a one-ended torsion-free hyperbolic group.

\begin{main}[Sela, Bowditch]\label{SB}
 Suppose that $Out(G)$ is an infinite group and that $G$ is not virtually a surface group. Then there is a non-trivial splitting of 
 $G$, given by a minimal action of $G$ on a tree $\mathscr{T}$ without 
 inversions and finite quotient $\mathscr{T}/G$. Moreover, the edge stabilizers are all cyclic subgroups and the set of vertex stabilizers 
 is $\aut(G)$-invariant. In particular, there is a finite index subgroup of $\aut(G)$ that preserves this splitting \textit{i.e.} 
 that preserves the conjugacy class of each vertex stabilizer and each edge stabilizer.
\end{main}

This splitting of $G$ is called the JSJ splitting.
The adjacency relation on the set of vertices induces a well-defined binary relation on the set of vertex stabilizers. 
The JSJ splitting is said canonical since $\aut(G)$ preserves the set of vertex stabilizers with the induced binary relation.
In fact, by considering the action of $G$ on $\partial G$, Bowditch proved in~\cite{Bowditch} that the extended action of $\aut(G)$
on $\mathscr{T}$ is $G$-compatible. An important ingredient of the work of Bowditch to prove Theorem \ref{SB} is due to 
Swarup~\cite{Swarup}. We refer to~\cite{Bowditch} and~\cite{Levitt} for a detailed description of the JSJ splitting. 
As in the case of a group that splits over a subgroup $C$, we obtain that the action of $G$ on $\mathscr{T}$ is minimal 
and without fixed points on the boundary.


\section{K\"ahler extensions and actions on trees}\label{AKG}

In this section we will prove our main results. 

\subsection{Applying Gromov and Schoen's Theorem}

Let $\Gamma$ be a K\"{a}hler group that fits in a short exact sequence 
\begin{gather}\label{eskh}
 \xymatrix{
 1\ar[r] & G \ar[r] &\Gamma \ar[r]^P &Q \ar[r] &1,
 }
\end{gather} 
where $G$ is a finitely generated group with trivial center acting on a tree $T$. Suppose that the action of $G$ on $T$ is minimal, 
faithful and without fixed points on the boundary. Note that in this case the action of $G$ on $T$ coincides with the induced action 
of $\inn(G)$ on $T$.

Let $\overline{\Gamma}$ be the image of the morphism  $\Gamma \rightarrow \aut(G)$ induced by the conjugation action of $\Gamma$ 
on $G$. The main result of this section is the following:

\begin{main}\label{kbs}
 Suppose that there is  a finite index subgroup $\overline{\Gamma}_0$ of $\overline{\Gamma}$ containing $\inn(G)$ such that
 the action of $G$ on $T$ can be extended to a $G$-compatible action of $\overline{\Gamma}_0$ on $T$. 
 Then $G$ is virtually a surface group. Moreover, there is a finite index subgroup $\Gamma_1$ of $\Gamma$ containing 
 $G$ such that the restricted short exact sequence 
  \begin{gather*}
 \xymatrix{
 1\ar[r] & G \ar[r] &\Gamma_1 \ar[r]^P &P(\Gamma_1) \ar[r] &1
 }
\end{gather*}
splits as a direct product.
\end{main}

To prove Theorem \ref{kbs} we will need a result of Gromov and Schoen~\cite{GS} on K\"{a}hler groups acting on trees. 
For the reader's convenience we recall some properties of holomorphic maps between compact complex manifolds and 
closed Riemann surfaces.
A surjective holomorphic map with connected fibers $f:X\rightarrow S$ from a compact complex manifold $X$ to a 
closed Riemann surface $S$ induces an orbifold structure as follows. For every point $p$ in $S$, let $m(p)$ be the greatest 
common divisor of the multiplicities of the irreducible components of the fiber $f^{-1}(p)$ and let $\Delta$ be the set of points 
in $S$ such that $m(p)\geq 2$. Note that $\Delta$ is finite since it is contained in the set of critical values of $f$. Hence, 
if $\Delta=\{p_1,...,p_k\}$ and $m_i=m(p_i)$ for $i=1,...,k$, the orbifold structure $\Sigma_f$ is given by 
$$\Sigma_f=\{S; (p_1,m_1),...,(p_k,m_k)\}.$$
For all $i=1,...,k$ let $\gamma_i$ be a loop around $p_i$ given by the boundary of a small enough disk such that $p_i$ is the only 
singular value contained in the disk. The fundamental group of the orbifold $\Sigma_f$ is given by 
$$\pi_1^{orb}(\Sigma_f)=\pi_1(S\setminus\Delta)//\ll\gamma_1^{m_1},...,\gamma_k^{m_k}\gg,$$
where $\ll\gamma_1^{m_1},...,\gamma_k^{m_k}\gg$ is the normal closure of $\{\gamma_1,...,\gamma_k\}$ in 
$\pi_1(S\setminus\Delta)$. 
The following lemma is a well-known result. A more general statement of this result can be found in~\cite{CKO}.

\begin{lemma}
 Let $X$ be a compact complex manifold, $S$ be a closed Riemann surface and $f:X\rightarrow S$ be a surjective holomorphic map whose 
 generic fiber $F$ is connected. If we denote by $\iota:F\hookrightarrow X$ the inclusion map and by $N$ the image of the induced 
 homomorphism on fundamental groups $\iota_*:\pi_1(F)\rightarrow \pi_1(X)$, we obtain a short exact sequence 
   \begin{gather*}
 \xymatrix{
 1\ar[r] & N \ar[r] &\pi_1(X) \ar[r]^{\Pi} &\pi_1^{orb}(\Sigma_f) \ar[r] &1.
 }
\end{gather*}
\end{lemma}

Since $F$ is compact we get that $N$ is finitely generated. The kernel of the map $\rho:\pi_1^{orb}(\Sigma_f)\rightarrow \pi_1(S)$ 
is isomorphic to $\ker(f_*)/N$ and we have the commutative diagram
\begin{gather*}
\xymatrix{
\pi_1(X)  \ar[rr]^{f_*} \ar[dr]_{\Pi} & & \pi_1(S)\raisetag{20mm} \\
& \pi_1^{orb}(\Sigma_f) \ar[ur]_{\rho}  \raisetag{20mm} \\
 }
\end{gather*}
To simplify the notation we will abbreviate $\Sigma_f$ to $\Sigma$. The result of Gromov and Schoen is the following.

\begin{main}[Gromov-Schoen]\label{t7}
 Let $X$ be a compact K\"{a}hler manifold whose fundamental group $\Gamma$ acts on a tree non-isomorphic to a 
 line nor a point. Suppose that the action is minimal with no fixed points on the boundary. Then there is a surjective
 holomorphic map with connected fibers from $X$ to a closed hyperbolic orbifold inducing the
 short exact sequence 
  \begin{gather*}
  \xymatrix{
  1\ar[r] &N \ar[r] &\Gamma \ar[r]^{\Pi\quad} &\pi_1^{orb}(\Sigma) \ar[r] &1,
  }
 \end{gather*}
 such that the restriction of the action to $N$ is trivial. 
\end{main}

{\noindent \it Proof of Theorem \ref{kbs}.} 
Let $\overline{\Gamma}_0$ be the finite index subgroup of $\overline{\Gamma}$ satisfying the hypothesis of the theorem and 
let $\Gamma_0$ be the preimage of $\overline{\Gamma}_0$ under the morphism $\Gamma \rightarrow \aut(G)$. Then, 
$ \Gamma_0$ is a finite index subgroup of $\Gamma$ containing $G$ such that the action of $G$ on $T$ can 
be extended to $\Gamma_0$.
 By Theorem \ref{t7} there is a short exact sequence
  \begin{gather*}
  \xymatrix{
  1\ar[r] &N \ar[r] &\Gamma_0 \ar[r]^{\Pi\quad} &\pi_1^{orb}(\Sigma) \ar[r] &1,
  }
 \end{gather*}
 where the restriction of the action of $\Gamma_0$ on $T$ to the subgroup $N$ is trivial and $\pi_1^{orb}(\Sigma)$ is virtually 
 isomorphic to a surface group.
 Now, let $\lambda:G\rightarrow\pi_1^{orb}(\Sigma)$ be the morphism given by the restriction of 
 $\Pi$ to $G$. This morphism is injective since $N\cap G$ is trivial. The latter assertion follows from the
 faithfulness of the action of $G$ on $T$. We claim that the subgroup $\Gamma_1=\Pi^{-1}(\lambda(G))$ 
 is the subgroup we are looking for. 
 First of all, since for a normal subgroup of $\pi_1^{orb}(\Sigma)$ being finitely generated is equivalent to having finite 
 index, we get that $\lambda(G)$ has finite index in $\pi_1^{orb}(\Sigma)$. This implies that
 $\Gamma_1$ is a finite index subgroup of $\Gamma_0$ (and thereby of $\Gamma$) and that $G$ is virtually isomorphic to a surface group. 
 To end the proof, it suffices to show that the morphism
 \begin{align*}
\eta: \Gamma_1\quad \rightarrow &\quad G\times P(\Gamma_1) \\
x \quad \mapsto &\, (\lambda^{-1}(\Pi(x)),P(x))
\end{align*}
 is bijective. The injectivity follows from the fact that the restriction of $\lambda^{-1}\circ\Pi$ to $G$ is the identity.
 For the surjectivity, if $(x_0,q_0)$ is an element of $G\times P(\Gamma_1)$, by taking 
 $y$ in $\Gamma_1$ such that $P(y)=q_0$ and
 $x$ in $G$ such that $\lambda(x)=\Pi(y)$, we get that $\eta(x_0x^{-1}y)=(x_0,q_0)$.
\hfill $\Box$

We used the fact that if the action of $G$ on $T$ is minimal without fixed points on the boundary, then any extension of this action 
has these properties as well. 


\subsection{Applications}\label{App}

Here we will apply Theorem \ref{kbs} and the results of Section \ref{GGT} to prove Theorems \ref{t}, \ref{kgfp}, \ref{kgbs} and 
\ref{kghg}. For this, we need to check that the action with which we start with is minimal, faithful and without fixed points on the 
boundary.
One can verify from the discussion at the end of Section \ref{BSTD}, that for a surface group (see also 
Section \ref{SG}), a free product of two groups and a Baumslag-Solitar group as in Theorem \ref{kgbs}, the actions 
on their respective Bass-Serre tree are minimal and without fixed points on the boundary of the tree.
For one-ended hyperbolic groups which are not virtually a surface group, 
this is a consequence of the theory of JSJ decompositions (see \cite{Bowditch,Sela}).
Hence, we only need to show faithfulness in all cases.

\begin{prop}\label{hg} 
 Let $S$ be a closed hyperbolic surface and let $C$ be the cyclic subgroup of $\pi_1(S)$ generated by the homotopy class of 
 a simple closed curve $\gamma$.
 Then the action of $\pi_1(S)$ on the Bass-Serre tree associated to the splitting of $\pi_1(S)$ over $C$ is faithful.
\end{prop}

\noindent {\it Proof.} 
Let $T$ be the Bass-Serre tree associated to the splitting of $\pi_1(S)$ over $C$.
Recall that $T$ coincides with the dual tree associated to $\gamma$ (see Remark \ref{dualtree}). This identification is given by a 
bijection between the edges of $T$ and a set of disjoint bi-infinite geodesics in $\hp$. Hence, an edge stabilizer of $T$ 
coincides with the stabilizer of a bi-infinite geodesic in $\hp$ under the action of $\pi_1(S)$  on $\hp$.
From this, we conclude that the intersection of the stabilizers of two different edges of $T$ is trivial, which implies the faithfulness 
of the action of $G$ on $T$.
\hfill $\Box$

{\noindent \it Proof of Theorem \ref{t}.} 
Let $\gamma$ be the simple closed curve in $S$, whose conjugacy class is preserved by the conjugation action of $\Gamma$ on 
$\pi_1(S)$. 
Let $C$ be the cyclic subgroup of $\pi_1(S)$ generated by the homotopy class of $\gamma$ and 
let $T$ be the Bass-Serre tree associated to the splitting of $\pi_1(S)$ over $C$. Therefore, the image of 
$\Gamma \rightarrow \aut(\pi_1(S))$  denoted by $\overline{\Gamma}$ is contained in the subgroup of automorphisms of $\pi_1(S)$ 
that preserves the conjugacy class of $\gamma$. By Proposition \ref{psg}, the action of $\pi_1(S)$ on $T$ can be extended to a 
$\pi_1(S)$-compatible action of a subgroup of $\overline{\Gamma}$ of index at most $2$ on $T$. 
By Proposition \ref{hg} we have that the action of $\pi_1(S)$ on $T$ is faithful and the result follows from Theorem \ref{kbs}.
\hfill $\Box$

Theorem \ref{kgfp} follows from Theorems \ref{fpt} and \ref{kbs}, and the following proposition (see also Remark \ref{fpfr}).

\begin{prop}\label{sgfp}
 If a surface group is embedded in a free product $A \ast B$, then it is embedded in $A$ or it is embedded in $B$ (up to conjugacy).
 In particular, a non-trivial free product of two groups is not virtually a surface group.
\end{prop}

\noindent {\it Proof.} 
 First, we will prove that a non-trivial free product is not isomorphic to a surface group. Let us assume by contradiction that 
 the free product $C \ast D$ is isomorphic to a surface group. Recall that an infinite index subgroup of a surface group is a free group. 
 Hence, since $C$ and $D$ are both infinite index subgroups of $C \ast D$, we obtain that they are both free groups. 
 Thus, the free product $C \ast D$ is also a free group, which is a contradiction. Note that the same argument shows that a 
 surface group is not isomorphic to a free product with more than two free factors.
 The general case follows from the Kurosh subgroup Theorem (see~\cite[p. 151]{SW}), which states that a subgroup of the free product 
 $A \ast B$ is given by the free product of a free group with subgroups of conjugates of $A$ and $B$. Hence, if a surface group 
 embeds in $A \ast B$, by the Kurosh subgroup Theorem and the previous argument it embeds as a subgroup of a conjugate of $A$ or as a 
 subgroup of a conjugate of $B$.
\hfill $\Box$
 
Theorem \ref{kgbs} follows from Theorem \ref{bsgaut} and \ref{kbs}, and the following two propositions:

\begin{prop}\label{pp1}
 The Baumslag-Solitar group 
 $$G(p,q)=\langle x,t | tx^{p}t^{-1}=x^{q}\rangle$$
 is not virtually a surface group.
\end{prop}

\noindent {\it Proof.} 
If $|p| \neq |q|$ and if we suppose that $G(p,q)$ is a surface group, we would have that $G(p,q)$ acts on $\hp$.
Hence, if we denote by $\ell$ the translation length of $x$ in $\hp$, as a consequence of the relation $tx^{p}t^{-1}=x^{q}$, we 
would have that $|p| \ell=|q| \ell$ which is impossible. 
If $H$ is a finite index subgroup of $G(p,q)$, there exists an integer $k$ such that $x^k$ and $t^k$ are contained in $H$. 
From the unique relation of the Baumslag-Solitar group, we deduce that $tx^{np}t^{-1}=x^{nq}$ for every integer $n$. Then, 
an argument by induction shows that 
$$t^{k}x^{k \cdot p^k}t^{-k}=x^{k \cdot q^k}.$$
Hence, if we suppose that $H$ is a surface group and we denote by $\ell$ the translation length of $x^k$ in $\hp$, we would have 
that $|p^k| \ell = |q^k| \ell$, which is impossible.
If $\vert p \vert = \vert q \vert$, the subgroup generated by $x^p$ is normal in $G(p,q)$. Hence its intersection with any finite 
index subgroup $H$ of $G(p,q)$ is normal in $H$. This cannot happen if $H$ is a surface group since such a group does not admit a 
non-trivial Abelian normal subgroup.
\hfill $\Box$

\begin{prop}\label{pp2}
 The Baumslag-Solitar group 
 $$G(p,q)=\langle x,t | tx^{p}t^{-1}=x^{q}\rangle,$$
 where $p, q$ are different integers with $p,q > 1$, acts faithfully on its Bass-Serre tree.
\end{prop}

\noindent {\it Proof.} 
The kernel $N$ of the action of $G(p,q)$ on its Bass-Serre tree is contained in every edge stabilizer. In particular, it is contained 
in $\langle x^{p}\rangle$. Assume by contradiction that $N$ is non-trivial. Then, $N$ is generated by $x^{kp}$ for some nonzero integer $k$.
Since $N$ is a normal subgroup of $G(p,q)$, we get that $tx^{kp}t^{-1}=x^{kq}$ is an element of $\langle x^{kp}\rangle$. 
Therefore $kq=mkp$ and $p$ divides $q$. By symmetry, we obtain that $q$ divides $p$ and thereby $p=q$, which is a contradiction.
\hfill $\Box$

Finally, Theorem \ref{kghg} follows from Theorems \ref{SB} and \ref{kbs}, and the fact that the action of a one-ended hyperbolic 
group with infinite group of outer automorphisms acts faithfully on its JSJ tree. The latter assertion follows from a 
similar argument to the one of Proposition \ref{hg} (in this case the intersection of the stabilizers of two different edges without 
common vertices is trivial).


\section{Surface groups and K\"{a}hler groups}\label{s4}


\subsection{K\"ahler extension of an Abelian group by a surface group}\label{qab}

Here we present a different proof of Bregman and Zhang's result using Theorem \ref{t}.

\begin{cor}[Bregman-Zhang]\label{c8}
 Let $S$ be a closed hyperbolic surface, $\Gamma$ be a K\"{a}hler group and $k$ be a positive integer such that there is a short exact 
 sequence 
  \begin{gather*}
  \xymatrix{
  1\ar[r] &\pi_1(S) \ar[r] &\Gamma \ar[r]^P &\Z^k \ar[r] &1.
  }
 \end{gather*}
Then there is a finite index subgroup $\Gamma_1$ of $\Gamma$ containing $\pi_1(S)$ such that the restricted short exact sequence
\begin{gather*}
  \xymatrix{
  1\ar[r] &\pi_1(S) \ar[r] &\Gamma_1 \ar[r]^{P\quad\quad} &P(\Gamma_1)\cong \Z^k \ar[r] &1
  }
 \end{gather*}
splits as a direct product.
 \end{cor}
 
 Note that $k$ must be even since $\Gamma_1$ is a K\"{a}hler group.
 Before proving Corollary \ref{c8}, we will recall some definitions and facts about the mapping class group of a closed hyperbolic 
 surface (this can be found with more details in~\cite{BLM}).
  Let $\mathscr{S}(S)$ denote the set of isotopy classes of simple closed curves in $S$. If $\tau$ is an element of the mapping 
 class group of $S$ and $\mathscr{A}$ is a subset of $\mathscr{S}(S)$, we denote by $\tau(\mathscr{A})$ the set 
 \[\{\tau(\alpha)|\alpha\in\mathscr{A}\},\]
 where $\tau(\alpha)$ denotes the isotopy class of $t(a)$ for $t\in\tau$ and $a\in\alpha$.
 A subset $\mathscr{A}$ of $\mathscr{S}(S)$ is said admissible if there is a set of disjoint simple closed curves that represent 
 all the isotopy classes of $\mathscr{A}$. An element $\tau$ of the mapping class group of $S$ is said to be 
 \begin{enumerate}
  \item {\it Reducible}. If there is a non-empty admissible set $\mathscr{A}$ such that $\tau(\mathscr{A})=\mathscr{A}$.
  \item {\it Pseudo-Anosov}. If for any isotopy class $\alpha$ in $\mathscr{S}(S)$ and for every nonzero integer $n$, 
  $\tau^n(\alpha)$ is different from $\alpha$ (this is one definition of a pseudo-Anosov mapping class among many others).
 \end{enumerate}
 To prove Corollary \ref{c8} we will need the classification theorem of Nielsen and Thurston and two lemmas.  
The first lemma is due to Birman, Lubotzky and McCarthy (see~\cite[Lemma~3.1]{BLM}).

\begin{lemma}[Birman-Lubotzky-McCarthy]\label{l10}
Let $A$ be an Abelian subgroup of the mapping class group of a hyperbolic surface generated by reducible elements 
$\{\tau_1,...,\tau_k\}$. Then, there is a non-empty canonical admissible set $\mathscr{A}$ such that 
$\tau_i(\mathscr{A})=\mathscr{A}$ for all $i=1,...,k$.
\end{lemma}
 
 The second lemma is the following. It will be a consequence of Thurston's hyperbolization Theorem (see~\cite{Otal}), 
 together with a classical result of Carlson and Toledo~\cite{CT}.
 
 \begin{lemma}\label{l11}
 Let $S$ be a hyperbolic surface, $\Gamma$ be a K\"{a}hler group and $k$ be a positive integer such that there is a short exact sequence
  \begin{gather*}
  \xymatrix{
  1\ar[r] &\pi_1(S) \ar[r] &\Gamma \ar[r]^{P} &\Z^k \ar[r] &1.
  }
 \end{gather*}
 Then the image of the monodromy map $\psi:\Z^k\rightarrow MCG(S)$ does not contain pseudo-Anosov elements.
 \end{lemma}
 
 \noindent {\it Proof.}  Let us assume by contradiction that the image of $\psi$ contains a pseudo-Anosov element $\tau$. 
 Since the elements of the centralizer of the cyclic subgroup generated by $\tau$ have the same set of fixed points on 
 the Thurston's boundary of the Teichm\"uller space of $S$, we get that 
 the cyclic subgroup generated by $\tau$ is a finite index subgroup of its centralizer (see~\cite[Theorem~4.6]{MP}, see 
 also~\cite{McCarthy}). Hence, by passing to a finite index subgroup of $\Gamma$ if necessary (which is also a K\"{a}hler 
 group), we may assume that the image of $\psi$ is the cyclic subgroup generated by $\tau$. Let $\{e_1,...,e_k\}$ be a basis 
 of $\Z^k$ such that $\{e_1,...,e_{k-1}\}$ generates the kernel of $\psi$ and $\psi(e_k)=\tau$. 
 We denote by $\Gamma_0$ the  subgroup of $\Gamma$ generated by $\pi_1(S)$ and $P^{-1}(e_k)$. 
 \begin{claim}\label{claimdp}
  $\Gamma$ is isomorphic to $\Gamma_0\times\Z^{k-1}$.
 \end{claim}
 Claim \ref{claimdp} implies the result since by Thurston's hyperbolization Theorem 
 (see~\cite{Otal}), $\Gamma_0$ is a cocompact lattice in the group of orientation preserving isometries of the hyperbolic $3$-space 
 and by Carlson and Toledo's Theorem (see~\cite{CT}), 
 the projection of $\Gamma$ onto the cocompact lattice $\Gamma_0$ factors through a surface group $\Lambda$. Hence, there is a 
 commutative diagram
\begin{gather*}
\xymatrix{
\Gamma\cong\Gamma_0\times\Z^k  \ar[rr]^{P_1} \ar[dr]_{\theta} & & \Gamma_0\raisetag{20mm} \\
& \Lambda \ar[ur]  \raisetag{20mm} \\
 } 
\end{gather*}
where $P_1$ is the projection onto $\Gamma_0$. This leads to the desired contradiction, since $\theta(\Gamma_0)$ is a subgroup 
of $\Lambda$ isomorphic to $\Gamma_0$ (as a consequence of the fact that $P_1|_{\Gamma_0}$ is the identity map), and every 
subgroup of $\Lambda$ is the fundamental group of a (closed or open) surface, but $\Gamma_0$ is neither free nor isomorphic to 
the fundamental group of a closed oriented surface since $H_3(\Gamma_0,\Z)$ is non-trivial.
\hfill $\Box$

\noindent {\it Proof of Claim \ref{claimdp}.}
Observe that for each $e_i$ in the kernel of $\psi$, there is a unique element $g_i$ in $\Gamma$ such that 
$g_i$ centralizes $\pi_1(S)$ and $P(g_i)=e_i$. Since the commutator $[g_i,g_j]$ is an element of $\pi_1(S)$ which centralizes 
$\pi_1(S)$ we get that the group $K$ generated by $\{g_1,...,g_{k-1}\}$ is isomorphic to $\Z^{k-1}$ and by construction 
$[K,\pi_1(S)]$ is trivial. Using the fact that $\pi_1(S)$ is a normal subgroup of $\Gamma$ and that $K$ centralizes $\pi_1(S)$ we 
obtain that $[K,\Gamma]$ centralizes $\pi_1(S)$. Finally, since $[K,\Gamma]$ is in the kernel of $P$ (which is $\pi_1(S)$) we conclude 
that $[K,\Gamma]$ is in the center of $\pi_1(S)$ which is trivial. Therefore the subgroups $K$ and $\Gamma_0$ commute and the result 
follows from observing that $\Gamma=K\Gamma_0$ and $K\cap\Gamma_0$ is trivial.
\hfill $\Box$

Now, we give our proof of Corollary \ref{c8}:

{\noindent \it Proof of Corollary \ref{c8}.} By Nielsen-Thurston's classification Theorem and by Lemma \ref{l11}, 
up to passing to a finite index subgroup of $\Z^k$ we have that the image of the monodromy map $\psi:\Z^k\rightarrow MCG(S)$ 
is generated by reducible elements. Hence, by Lemma \ref{l10} there is a non-empty canonical admissible set $\mathscr{A}$ of 
$\mathscr{S}(S)$ such that $\psi(\Z^k)$ preserves $\mathscr{A}$. Since a finite index subgroup of $\Z^k$ (which is isomorphic 
to $\Z^k$) acts trivially on $\mathscr{A}$, we can apply Theorem \ref{t} which concludes the proof.
\hfill $\Box$


\subsection{K\"ahler extension by a surface group whose monodromy map has Abelian image}
With a bit more work we get the next corollary of Theorem \ref{t} which is more general than Corollary \ref{c8}.

\begin{cor}\label{c12}
Let $S$ be a hyperbolic surface and $\Gamma$ be a K\"{a}hler group such that there is a short exact sequence 
  \begin{gather*}
  \xymatrix{
  1\ar[r] &\pi_1(S) \ar[r] &\Gamma \ar[r]^{P} &Q \ar[r] &1.
  }
 \end{gather*}
 If the image of the monodromy map $\psi:Q\rightarrow MCG(S)$ is Abelian, then there are finite index subgroups $\Gamma_1$ of $\Gamma$ 
 and $Q_1$ of $Q$ such that $\pi_1(S)$ is contained in $\Gamma_1$ and the restricted short exact sequence
 \begin{gather*}
  \xymatrix{
  1\ar[r] &\pi_1(S) \ar[r] &\Gamma_1 \ar[r] &Q_1 \ar[r] &1
  }
 \end{gather*}
 splits as a direct product.
\end{cor}

To prove Corollary \ref{c12} it suffices to prove a result analogous to Lemma \ref{l11} for the case when
the monodromy map $\psi:Q\rightarrow MCG(S)$ has Abelian image (instead of assuming that Q itself is Abelian).
According to a result of Birman, Lubotzky and McCarthy (see~\cite[Theorem~B]{BLM}),every solvable subgroup of the mapping
class group of a hyperbolic surface is virtually Abelian. Hence, 
Corollary \ref{c12} can be extended to the case when $\psi(Q)$ is a solvable group.

\begin{lemma}\label{l14}
Let $S$ be a hyperbolic surface and $\Gamma$ be a K\"{a}hler group such that there is a short exact sequence
  \begin{gather*}
  \xymatrix{
  1\ar[r] &\pi_1(S) \ar[r] &\Gamma \ar[r]^{P} &Q \ar[r] &1.
  }
 \end{gather*}
 If the image of the monodromy map $\psi:Q\rightarrow MCG(S)$ is Abelian, then it does not contain pseudo-Anosov
  elements.
\end{lemma}
\noindent {\it Proof.}  
The proof is by contradiction. Let us suppose that the image of $\psi$ contains a pseudo-Anosov element $\tau$. As in 
Lemma \ref{l11}, up to passing to a finite index subgroup of $\Gamma$ we may assume that the image of $\psi$ is the cyclic 
subgroup generated by $\tau$. Let $t$ be an element in $\Gamma$ such that $\psi\circ P(t)=\tau$ and let $\Gamma_0$ be the  subgroup 
of $\Gamma$ generated by $\pi_1(S)$ and $t$. The result will follow from exhibiting  a homomorphism $\Gamma\rightarrow\Gamma_0$ 
whose restriction to $\Gamma_0$ is the identity, since Thurston's hyperbolization Theorem  and 
Carlson and Toledo's Theorem will lead us to a contradiction as in Lemma \ref{l11}.

Now, to construct this homomorphism, observe that since each element in the kernel of $\psi$ has a unique lift in $\Gamma$ that 
centralizes $\pi_1(S)$, there is a subgroup $K$ of $\Gamma$ isomorphic to the kernel of $\psi$ that centralizes $\pi_1(S)$. 
Indeed, $K$ is the centralizer of $\pi_1(S)$ in $\Gamma$. Hence, 
the subgroup $[\Gamma,K]$ centralizes $\pi_1(S)$ and since its image under $P$ is contained in the kernel of $\psi$, we get that 
$[\Gamma,K]$ is contained in $K$ and therefore $K$ is a normal subgroup of $\Gamma$. Finally, since the kernel 
of $\psi\circ P$ is isomorphic to $K\times\pi_1(S)$,
we obtain that $\Gamma$ is isomorphic to the semi-direct product $(K\times\pi_1(S))\rtimes\langle t \rangle$, and 
thereby there is a well-defined homomorphism 
 \begin{align*}
\Gamma \rightarrow &\Gamma_0 \\
kxt^m \mapsto &\, xt^m,
\end{align*}
where $k$ is in $K$, $x$ is in $\pi_1(S)$ and $m$ is in $\Z$. To verify that it is a well-defined map we just need to observe that 
for $k_1x_1t^{m_1}$ and $k_2x_2t^{m_2}$ in $\Gamma$ we have that 
\[k_1x_1t^{m_1}k_2x_2t^{m_2}=(k_1t^{m_1}k_2t^{-m_1})(x_1t^{m_1}x_2t^{-m_1})t^{m_1+m_2},\]
where we notice that $t^{m_1}k_2t^{-m_1}$ is in $K$ and $t^{m_1}x_2t^{-m_1}$ is in $\pi_1(S)$.
\hfill $\Box$

The following result is a consequence of a theorem due to Ivanov which states the following. 
Let $S$ be  a closed hyperbolic surface and $\Gamma$  bean infinite subgroup of $MCG(S)$ that 
does not preserve any admissible set $\mathscr{A}$ of $\mathscr{S}(S)$. Then, $\Gamma$ is
either virtually a cyclic group generated by a pseudo-Anosov element, or 
it contains a free group generated by $2$ pseudo-Anosov elements (see~\cite[Theorem~2]{Ivanov}).

\begin{cor}\label{l19}
Let $S$ be a hyperbolic surface and $\Gamma$ be a K\"{a}hler group such that there is a short exact sequence
  \begin{gather*}
  \xymatrix{
  1\ar[r] &\pi_1(S) \ar[r] &\Gamma \ar[r]^{P} &Q \ar[r] &1.
  }
 \end{gather*}
 If the image of the monodromy map $\psi:Q\rightarrow MCG(S)$ is infinite, then it contains a free subgroup 
 generated by $2$ pseudo-Anosov elements.
\end{cor}
\noindent {\it Proof.} 
Recall that an extension as above is virtually trivial if and only if the monodromy 
subgroup of the short exact sequence is finite. Therefore, by Theorem \ref{t}, $\psi(Q)$ 
does not preserve any admissible set $\mathscr{A}$ of $\mathscr{S}(S)$.  Finally, the result follows from Lemma \ref{l14}
and Ivanov's result recalled before the statement.
\hfill$\Box$


\subsection{More restrictions on K\"ahler extensions by a surface group}

A classical result due to Scott (see~\cite{Scott}) allows us to extend Theorem \ref{t} to any closed curve in $S$. 

\begin{main}[Scott]\label{t15}
Let $S$ be a topological surface and let $x$ be an element of $\pi_1(S)$. Then, there is a finite covering map $S'\rightarrow S$
such that $x$ is in $\pi_1(S')$ and can be represented by a simple closed curve in $S'$.
\end{main}

\begin{lemma}\label{l17}
Let $G$ be a group and $H$ be a finitely generated normal subgroup of $G$. If $\Lambda$ is a finite index subgroup 
of $H$, then the normalizer of $\Lambda$ in $G$ is a finite index subgroup of $G$.
\end{lemma}
\noindent {\it Proof.}  Let $n$ be the index of $\Lambda$ in $H$. Then $G$ acts on the set of subgroups of $H$ of index 
$n$ by conjugation. This set is finite and therefore the stabilizer of any of these subgroups 
gives a finite index subgroup of $G$. Finally, notice that the stabilizer of $\Lambda$ is precisely the normalizer of 
$\Lambda$ in $G$ and the result follows.
\hfill$\Box$

\begin{main}\label{t18}
 Let $\Gamma$ be a K\"{a}hler group such that there is a short exact sequence 
 \begin{gather*}
  \xymatrix{
  1\ar[r] &\pi_1(S) \ar[r] &\Gamma \ar[r]^{P} &Q \ar[r] &1,
  }
 \end{gather*}
 with $S$ a closed hyperbolic surface. If the conjugation action of $\Gamma$ on $\pi_1(S)$ preserves the conjugacy class of a 
 non-trivial element of  $\pi_1(S)$, then there is a finite index subgroup $\Gamma_1$ of $\Gamma$ containing a finite index subgroup 
 $\Lambda$ of $\pi_1(S)$ which is normal in $\Gamma_1$ such that the extension 
 \begin{gather*}
  \xymatrix{
  1\ar[r] &\Lambda \ar[r] &\Gamma_1 \ar[r] &\Gamma_1/\Lambda \ar[r] &1
  }
 \end{gather*}
 splits as a direct product.
\end{main}
\noindent {\it Proof.}  
Let $\gamma$ be the closed curve in $S$ whose conjugacy class is preserved by the action of $\Gamma$. 
We may assume that $\gamma$ is not simple. By Theorem \ref{t15}, there is a finite covering map $S'\rightarrow S$, such that 
$\gamma$ lifts to a simple closed curve in $S'$. 
Let $\Lambda$ be the fundamental group of $S'$
and let $\Gamma'$ be the normalizer of $\Lambda$ in $\Gamma$.
By Lemma \ref{l17}, $\Gamma'$ is a finite index 
subgroup of $\Gamma$ and therefore $\Gamma'$ is a K\"{a}hler group. Hence, the short exact sequence 
 \begin{gather*}
  \xymatrix{
  1\ar[r] &\Lambda \ar[r] &\Gamma' \ar[r] &\Gamma'/\Lambda \ar[r] &1
  }
 \end{gather*}
 satisfies the hypothesis of Theorem \ref{t}, which implies the existence of a subgroup $\Gamma_1$ of $\Gamma'$ of finite index (and 
 therefore a finite index subgroup of $\Gamma$) such that the exact sequence 
  \begin{gather*}
  \xymatrix{
  1\ar[r] &\Lambda \ar[r] &\Gamma_1 \ar[r] &\Gamma_1/\Lambda \ar[r] &1,
  }
 \end{gather*}
 splits as a direct product.
 \hfill$\Box$


{\bf Acknowledgements.} This article is a part of my PhD thesis. I would like to thank my PhD advisors, Thomas Delzant and Pierre Py who 
proposed this subject to me, for all the discussions concerning this work and for their great support. 
I would also like to thank Beno\^it Claudon, Claudio Llosa Isenrich and the referee for their comments on this text, and 
Martin Mion-Mouton and Ad\`ele P\'erus for reading a preliminary version of this text.

\bigskip
\bigskip
\begin{small}
\begin{tabular}{lllll}
Francisco Nicol\'as Cardona & & & &\\
IRMA & & & & \\
Universit\'e de Strasbourg \& CNRS & & & & \\
67084 Strasbourg, France & & & & \\
fnicolascardona@math.unistra.fr & & & &\\    
\end{tabular}
\end{small}


\begin{thebibliography}{00}

\bibitem{ABCKT} J.~Amor\'{o}s, M.~Burger, K.~Corlette, D.~Kotschick and D.~Toledo, \emph{Fundamental groups of compact 
{K}\"{a}hler manifolds}, Mathematical Surveys and Monographs.~{\bf 44}, American Mathematical Society, Providence, RI, (1996).

\bibitem{BLM} J. S.~Birman,  A.~Lubotzky and J.~McCarthy,  \emph{Abelian and solvable subgroups of the mapping class groups},
Duke Mathematical Journal.~{\bf 50}, No.~4, (1983), 1107--1120.

\bibitem{Bowditch} B. H.~Bowditch, \emph{Cut points and canonical splittings of hyperbolic groups}, Acta Mathematica.~ {\bf 180}, 
No.~2, (1998), 145--186.

\bibitem{BZ} C.~Bregman and L.~Zhang, \emph{On K\"{a}hler extensions of abelian groups}, arXiv e-prints, (2016).

\bibitem{CT} J. A.~Carlson and D.~Toledo, \emph{Harmonic mappings of {K}\"{a}hler manifolds to locally symmetric spaces},
Institut des Hautes \'{E}tudes Scientifiques. Publications Math\'{e}matiques, No.~69, (1989), 173--201.

\bibitem{CKO} F.~Catanese, J.~Keum and K.~Oguiso, \emph{Some remarks on the universal cover of an open {$K3$} surface},
Mathematische Annalen, {\bf 325}, No.~2, (2003), 279--286.

\bibitem{Dehn} M.~Dehn, \emph{Papers on group theory and topology}, Translated from the German and with introductions and an
appendix by John Stillwell, With an appendix by Otto Schreier, Springer-Verlag, New York, (1987).


\bibitem{FR1} D. I.~Fouxe-Rabinovitch, \emph{\"{U}ber die {A}utomorphismengruppen der freien {P}rodukte. {I}},
Rec. Math. [Mat. Sbornik] N.S., {\bf 8 (50)}, (1940), 265--276.

\bibitem{FR2} D. I.~Fouxe-Rabinovitch, \emph{\"{U}ber die {A}utomorphismengruppen der freien {P}rodukte. {II}},
Rec. Math. [Mat. Sbornik] N. S., {\bf 9 (51)}, (1941), 183--220.

\bibitem{GHMR} N. D.~Gilbert, J.~Howie, V.~Metaftsis and E.~Raptis, \emph{Tree actions of automorphism groups},
Journal of Group Theory.~{\bf 3}, No.~2, (2000), 213--223.

\bibitem{GS} M.~Gromov and R.~Schoen, \emph{Harmonic maps into singular spaces and {$p$}-adic superrigidity for lattices in 
groups of rank one}, Institut des Hautes \'{E}tudes Scientifiques. Publications Math\'{e}matiques, No.~76, (1992), 165--246.

\bibitem {Ivanov} N. V.~Ivanov, \emph{Subgroups of {T}eichm\"{u}ller modular groups}, Translations of Mathematical Monographs.~
{\bf 115}, Translated from the Russian by E. J. F. Primrose and revised by the author, American Mathematical Society, Providence, RI,
(1992).

\bibitem{KPS} A.~Karrass, A.~Pietrowski and D.~Solitar, \emph{Automorphisms of a free product with an amalgamated subgroup},
Contributions to group theory, Contemp. Math., Amer. Math. Soc., Providence, RI.~{\bf 33}, (1984), 328--340.

\bibitem{Levitt} G.~Levitt, \emph{Automorphisms of hyperbolic groups and graphs of groups}, Geometriae Dedicata.~ {\bf 114},
(2005), 49--70.


\bibitem{McCarthy} J.~McCarthy, \emph{A ``{T}its-alternative'' for subgroups of surface mapping class groups},
Transactions of the American Mathematical Society.~{\bf 291}, No.~ 2, (1985), 583--612.

\bibitem{MP} J.~McCarthy and A.~Papadopoulos, \emph{Dynamics on {T}hurston's sphere of projective measured foliations},
Commentarii Mathematici Helvetici.~{\bf 64}, No.~1, (1989), 133--166.

\bibitem{Nielsen} J.~Nielsen, \emph{Untersuchungen zur {T}opologie der geschlossenen zweiseitigen {F}l\"{a}chen}, 
Acta Mathematica.~{\bf 50}, No.~1, (1927), 189--358.

\bibitem{Otal} J. P.~Otal, \emph{Le th\'{e}or\`eme d'hyperbolisation pour les vari\'{e}t\'{e}s fibr\'{e}es de dimension 3},
Ast\'{e}risque, No.~235, (1996).

\bibitem{Pettet} M. R.~Pettet, \emph{The automorphism group of a graph product of groups}, Communications in Algebra.~{\bf 27}, 
No.~10, (1999), 691--4708.

\bibitem{Scott} P.~Scott, \emph{Subgroups of surface groups are almost geometric}, Journal of the London Mathematical Society. 
Second Series {\bf 17}, No.~3, (1978), 555--565.

\bibitem{SW} P.~Scott and T.~Wall, \emph{Topological methods in group theory}, Homological group theory ({P}roc. {S}ympos., 
{D}urham, 1977), London Math. Soc. Lecture Note Ser., Cambridge Univ. Press, Cambridge-New York.~{\bf 36}, (1979), 137--203.

\bibitem{Sela} Z.~Sela, \emph{Structure and rigidity in ({G}romov) hyperbolic groups and discrete groups in rank {$1$} {L}ie groups. {II}},
Geometric and Functional Analysis.~ {\bf 7}, No.~3, (1997), 561--593.

\bibitem{Serre} J. P.~Serre, \emph{Trees}, Springer Monographs in Mathematics, Translated from the French original by John Stillwell,
Corrected 2nd printing of the 1980 English translation, Springer-Verlag, Berlin, (2003).
 
 \bibitem{Shiga} H.~Shiga, \emph{On monodromies of holomorphic families of {R}iemann surfaces and modular transformations},
Mathematical Proceedings of the Cambridge Philosophical Society.~{\bf 122}, No.~3, (1997), 541--549.

\bibitem{Sun} X.~Sun, \emph{Regularity of harmonic maps to trees}, American Journal of Mathematics.~{\bf 125}, No.~4, (2003), 737--771.

\bibitem{Swarup} G. A.~Swarup, \emph{On the cut point conjecture}, Electronic Research Announcements of the American Mathematical
Society.~{\bf 2}, No.~2, (1996), 98--100.


\end{thebibliography}
\end{document}